\documentclass[conference]{IEEEtran}
\usepackage[]{times}
\usepackage{epsfig}
\usepackage{subfigure}
\usepackage{amsmath}
\usepackage{bm,amssymb}
\usepackage{graphicx}
\usepackage[top=0.5in, bottom=0.75in, left=0.75in, right=0.75in]{geometry}
\usepackage{multirow}
\usepackage{flushend}
\usepackage{url}
\usepackage{color}

\title{Distributed control of reactive power flow\\ in a radial distribution circuit with high photovoltaic penetration}
\author{
\authorblockN{Konstantin Turitsyn}
\authorblockA{CNLS \& Theoretical Divison\\
Los Alamos National Lab\\
Los Alamos,\\ NM 87545, USA\\
Email: turitsyn@lanl.gov}
\and
\authorblockN{Petr \v{S}ulc}
\authorblockA{New Mexico Consortium,\\ Los Alamos,\\ NM
87544, USA\\Email: sulcpetr@gmail.com}
\and
\authorblockN{Scott Backhaus}
\authorblockA{Materials, Physics\\ \& Applications Division\\
Los Alamos National Lab\\
Los Alamos,\\ NM 87545, USA\\
Email: backhaus@lanl.gov}
\and
\authorblockN{Michael Chertkov}
\authorblockA{CNLS \& Theoretical Divison\\
Los Alamos National Lab\\
Los Alamos,\\ NM 87545, USA\\
Also with NMC\\
Email: chertkov@lanl.gov}}

\begin{document}
\maketitle
\begin{abstract}
We show how distributed control of reactive power can serve to regulate voltage and minimize resistive losses in a distribution circuit that includes a significant level of photovoltaic (PV) generation.   To demonstrate the technique, we consider a radial distribution circuit with a single branch consisting of sequentially-arranged residential-scale loads that consume both real and reactive power.  In parallel, some loads also have PV generation capability. We postulate that the inverters associated with each PV system are also capable of limited reactive power generation or consumption, and we seek to find the optimal dispatch of each inverter's reactive power to both maintain the voltage within an acceptable range and minimize the resistive losses over the entire circuit.  We assume the complex impedance of the distribution circuit links and the instantaneous load and PV generation at each load are known.  We compare the results of the optimal dispatch with a suboptimal local scheme that does not require any communication.  On our model distribution circuit, we illustrate the feasibility of
 high levels of PV penetration and a significant (20\% or higher) reduction in losses.

\end{abstract}


\section{Introduction}
\label{sec:intro}

Utilities employ various equipment to control the operation of primary distribution systems including under-load tap changing transformers (ULTC), step voltage regulators (SVR), and fixed and switchable capacitors (FC and SC).  The primary function of these devices is to maintain voltage at the customer service entrance within an acceptable range to ensure adequate operation and lifetime of customer equipment.  In many distribution systems, the operation of these devices is governed by local conditions primarily local voltage and current sensors.  Utilities are increasingly utilizing distribution-level Supervisory Control and Data Acquisition (SCADA) to operate this equipment, and in some cases, this has led to increased levels of centralized control of this equipment allowing for coordinated operation.  Centralized, coordinated control also provides the opportunity to optimize the operation to meet utility goals such as minimization of losses, reduction of peak apparent power, or extension of equipment life.  A review of previous work in optimal placement and sizing of FCs and SCs and coordinated operation of SCs and ULTCs can be found series of papers by Baran and Wu~\cite{89BWb,89BWa} and Baldick and Wu~\cite{90BW}.

The relatively slow operation of ULTCs, SVRs, and SCs is acceptable in most distribution systems where fluctuations in loads and voltage levels are relatively small and significant changes in average load occur relatively slowly and in a predictable fashion through out the day and year.  However, this situation is due to change as distribution systems are subject to higher levels of time-variable distributed renewable generation, primarily residential solar photovoltaic (PV) systems.  High penetrations of time-variable renewable generation will pose several new challenges to voltage regulation but may also create new opportunities for optimization of distribution systems.

The distributed nature of residential PV generation implies that real power is injected at many points along the circuit making it no longer possible to obtain reliable estimates of the power flows throughout the circuit from a few measurements of current made at a few discrete locations.  In addition, real power flows will increasingly be two way and in some circumstances, the present deployment of ULTCs, SVRs, and SCs, may not be sufficient to ensure adequate voltage regulation at the service entrances.\cite{08BMOCKS}  Furthermore, the variability of PV generation can occur on a timescale much shorter than the present equipment can cope with.  For instance, cloud transients can cause ramps in PV generation on the order of 15\% per second at a particular location slowing to perhaps 15\% per minute for an entire distribution circuit due to its spatial diversity.\cite{90GKFE,08WNRN} At high levels of PV penetration on a distribution circuit, this can result in a reversal of real power flow (i.e. change from net generation to net consumption) over a period of a few minutes and a loss of voltage regulation due to the slow response of existing equipment. The existing equipment could respond on shorter timescales, however, the increased number of operations that would be required to counteract the variability due to weather conditions would drastically reduce the lifetime of the switches and tap changers \cite{08BMOCKS,08WNRN}.

To mitigate many of the issues discussed above, it has been proposed that the interconnection standards for inverter-based distributed generation be changed in such a way to enable the inverters to assist with high speed voltage regulation\cite{08BMOCKS,08WNRN,08LB}.  To be compliant with present interconnection standards,\cite{1547} PV inverters must not inject or consume reactive power or in anyway attempt to regulate voltage.  However, inverters with this control capability already exist for off-grid applications and for grid-tied applications when the PV system is operating in an islanded mode.  Availability of hardware is not a fundamental barrier, but excess apparent power capacity (above the real power capability of the PV system it is connected to) must be built into an inverter to allow for reactive power generation and consumption while operating near maximum real power.  Determining the appropriate size of this additional capacity is an important outstanding question
and depends on yet-to-be-developed control schemes that will coordinate the inverters' response to changes in voltage and power flow.  In addition to effectively regulating voltage, a coordinated control scheme should also allow for the optimization of the reactive power flows to minimize dissipation in the distribution circuit. This paper presents a case study of such an algorithm utilizing optimal and distributed control of PV-inverter reactive power generation.

The layout of the material in the remainder of this manuscript is as follows. Section
\ref{sec:Inverter} describes a simplified model of an inverter capable of limited reactive power generation and consumption. Section \ref{sec:DistFlow} describes an optimization problem where we utilize the inverters' additional capacity to minimize losses in a radial distribution circuit while respecting the constraints of voltage regulation and the inverters' apparent power capacity.  In Section \ref{sec:rural}, we describe the parameters for a model distribution circuit that serves as a prototype of a sparsely-loaded rural distribution circuit. Section \ref{sec:simulations} reports the simulation results demonstrating feasibility of the distributed control and illustrating the quality of improvements possible on our prototype circuit. Finally, Section
\ref{sec:conclusions} discusses our conclusions and path forward.

\section{Inverter as a limited regulator of local reactive power flow.}
\label{sec:Inverter}

\begin{figure}
\centering \caption{When $s$ is larger than $p^{(g)}$, the inverter can supply or consume reactive power $q^{(g)}$.  The inverter can dispatch $q^{(g)}$ quickly (on the cycle-to-cycle time scale) providing a mechanism for rapid voltage regulation.  As the output of the PV panel array $p^{(g)}$ approaches $s$, the range of available $q^{(g)}$ decreases to zero.
} \includegraphics[width=0.5\textwidth]{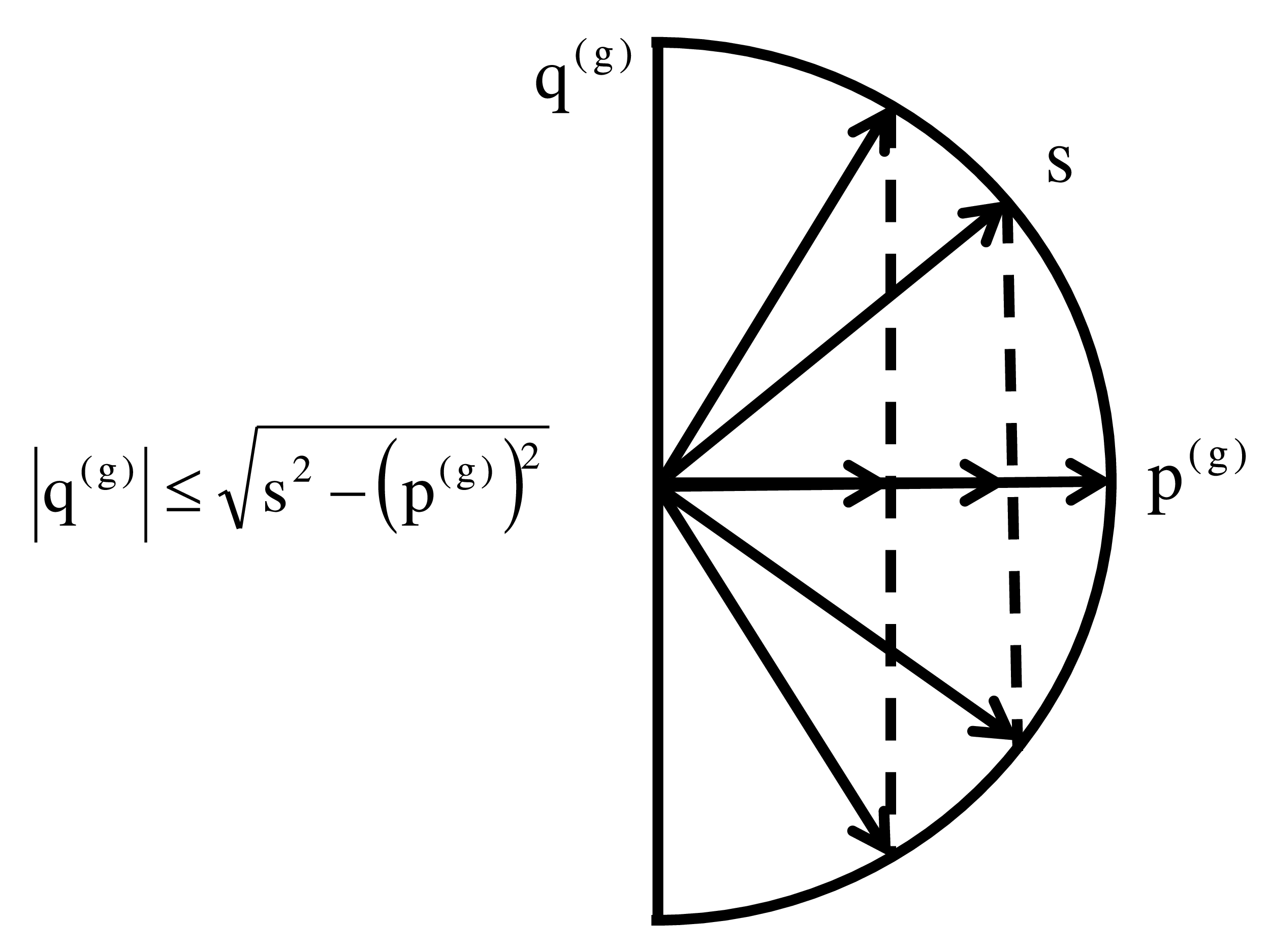}
\label{fig:complex}
\end{figure}

The present interconnection standard for PV-inverters \cite{1547} forces all inverters to operate at unity power factor while operating in a grid-tied mode, i.e. the inverter must not generate or consume reactive power or attempt to regulate voltage.   Several researchers have proposed that, to allow for high penetrations of PV on a distribution circuit, the present standard is not workable because the "utility-scale" regulation equipment discussed earlier is not sufficient to handle two-way power flows or fast enough to mitigate rapid cloud transients \cite{08BMOCKS,08WNRN,08LB}.  To be used in this capacity, PV inverters must be allowed a new degree of freedom to provide voltage regulation.  A possible scenario is that a PV inverter will be required to have a maximum apparent power capability $s$ larger than the maximum power output of its PV panel array, $\max p^{(g)}$, and the excess capability will be dispatched by the distribution utility to provide for voltage regulation.

A simple model for the relationship between the various PV inverter output variables has been described previously.\cite{08LB}  If $s$ is larger than $p^{(g)}$, the inverter can supply or consume reactive power $q^{(g)}$.  The magnitude of $q^{(g)}$ is bounded by $\sqrt{s^2-(p^{(g)})^2}$ and decreases as the real power output of the inverter approaches $s$.  The phasor relationship between the inverter operating parameters is shown in Fig. \ref{fig:complex} for several different levels of PV panel power output, $p^{(g)}$.

Although advanced inverters may have the capability to generate and consume reactive power, this output must be dispatched in such a way to effectively regulate voltage and achieve other utility goals.  The algorithms used to perform the dispatch could be based on either local conditions or be done centrally using circuit-wide information.  Because of the diversity of distribution utility infrastructure, operating areas, and economic models, we develop and compare both approaches.

To deploy centralized algorithms, we envision that the communication and control functions of smart-grid technologies will be crucial.  Enhanced communications capabilities of advanced smart-grid systems will allow distribution automation systems (DAS) high speed access to voltage amplitude and power measurements at many if not all service entrances greatly expanding the DAS's grid visibility and situational awareness.  These measurements, coupled with power flow models of the circuit, will allow optimization algorithms (developed here and by others) running in the DAS to individually dispatch reactive power from each PV inverter to ensure that service voltages stay within acceptable bounds and, for instance, to minimize losses.  With this level of system knowledge and centralized approach, algorithms can be developed that guarantee optimal solutions.

However, not all smart-grid schemes allow for high speed communication, e.g. power line carrier (PLC) schemes may be too slow to update voltage and power measurements fast enough to mitigate the effects of cloud transients.  In this case, local schemes that only rely upon local measurements, perhaps including knowledge of nearest-neighbor nodes, are developed.  In this article, we develop both an optimal centralized and suboptimal but local algorithm and compare their performance.  The results demonstrate that the performance of a simple distributed scheme can approach the performance of centralized schemes.

\section{Power Flow. Optimization of Losses and Voltage Control.}
\label{sec:DistFlow}

\begin{figure}
\centering \caption{
Diagram and notations for the radial network. $P_j$ and $Q_j$ represent real and reactive power flowing down the circuit from node $j$, where $P_0$ and $Q_0$ represent the power flow from the sub-station. $p_j$ and $q_j$ correspond to the flow of power out of the network at the node $j$, where the respective positive [negative] contributions, $p_j^{(c)}$ and $q_j^{(c)}$ [$p_j^{(g)}$ and $q_j^{(g)}$] represent consumption [generation] of power at the node. The node-local control parameter $q_j^{(p)}$ can be positive or negative but is bounded in absolute value as described in Eq.~\ref{PV_constraint}. The apparent power capability of the inverter $s_j$ is preset to a value comparable to but larger than $\max p_j^{(g)}$}. \includegraphics[width=0.5\textwidth]{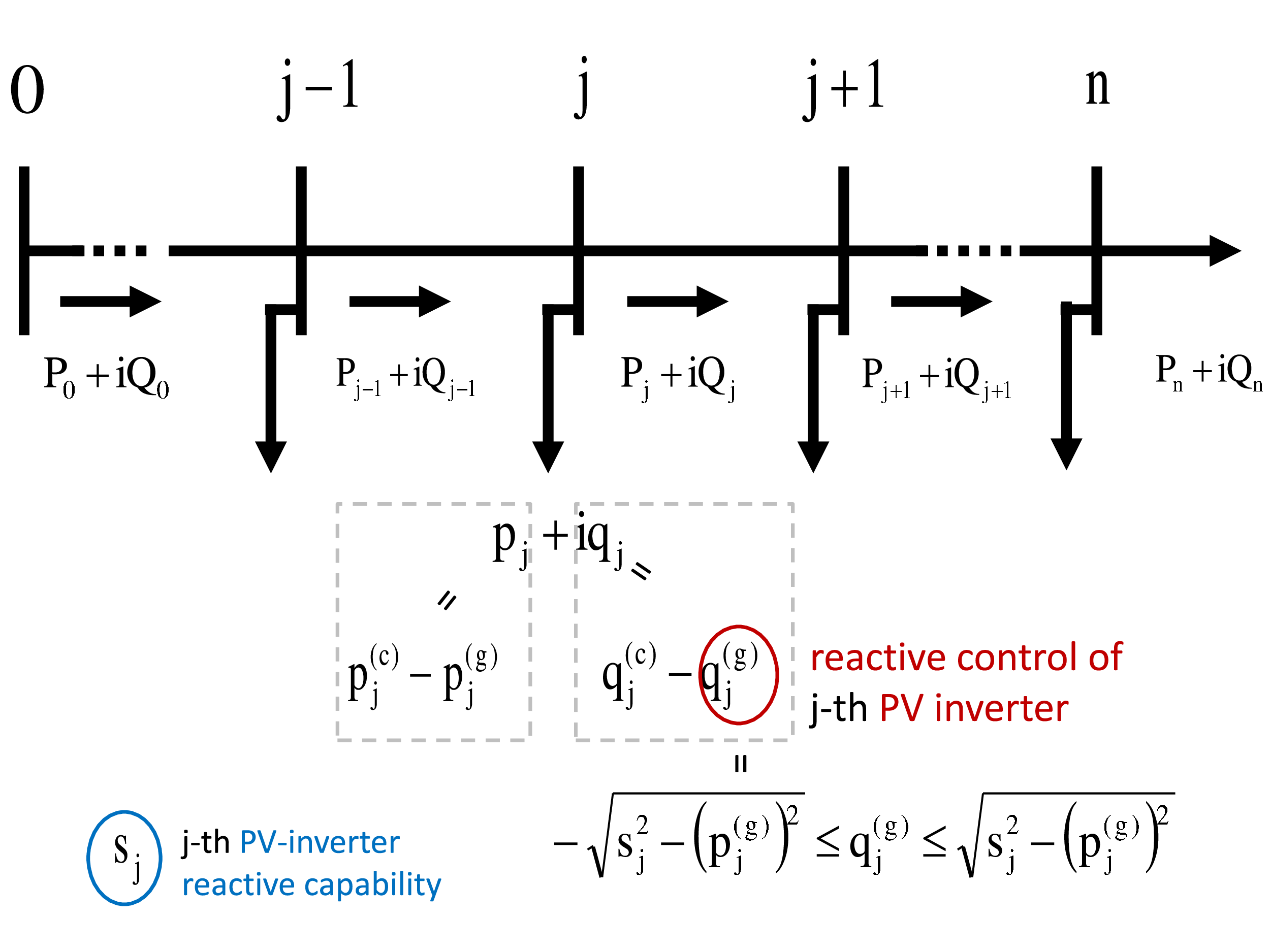}
\label{fig:feeder}
\end{figure}

To solve for complex power flows in a radial circuit, we follow  the DistFlow description of \cite{89BWb,89BWa,89BWc}. The system of AC power flow equations representing the radial circuit illustrated in Fig.~(\ref{fig:feeder}) is $\forall j=1,\cdots,n$:
\begin{eqnarray}
&&P_{j+1}\!=\! P_j\!-\!r_j\frac{P_j^2\!+\!Q_j^2}{V_j^2}\!-\!p_{j+1}, \label{Pj+1}\\
&&Q_{j+1}\!=\!Q_j\!-\!x_j\frac{P_j^2\!+\!Q_j^2}{V_j^2}\!-\!q_{j+1}, \label{Qj+1}\\
&&V_{j+1}^2\!=\!V_j^2\!-\!2(r_jP_j\!+\!x_jQ_j)\!+\!(r_j^2\!+\!x_j^2)
\frac{P_j^2\!+\!Q_j^2}{V_j^2},
\label{Vj2}
\end{eqnarray}
where $P_j+iQ_j$ is the complex power flowing away from node $j$ toward node $j+1$, $V_j$ is the voltage at node $j$, $r_j+ix_j$ is the complex impedance of the link between node $j$ and $j+1$, and $p_j+i q_j$ is the complex power extracted at the node $j$. Both $p_j$ and $q_j$ are composed of local consumption minus local generation due to the PV inverter, i.e. $p_j=p_j^{(c)}-p_j^{(g)}$ and $q_j=q_j^{(c)}-q_j^{(g)}$.  Of the four contributions to $p_j+i q_j$, $p_j^{(c)}$, $p_j^{(g)}$ and, $q_j^{(c)}$ are uncontrolled (i.e. driven by consumer load or instantaneous PV generation), while the reactive power generated by the PV inverter, $q_j^{(g)}$, can be adjusted.   However, $q_j^{(g)}$ is limited by the reactive capability of the inverter:
\begin{eqnarray}
\forall j=1,\cdots,n:\quad
\left| q_j^{(g)} \right| \leq \sqrt{s_j^2-(p_j^{(g)})^2}. \label{PV_constraint}
\end{eqnarray}

The rate of energy dissipation (losses) in the distribution circuit,
\begin{eqnarray}
{\cal L}=\sum_{j=0}^{n-1} r_j\frac{P_j^2+Q_j^2}{V_j^2},
\label{Loss}
\end{eqnarray}
is an important global characteristic.  Minimizing or at least keeping the losses acceptably low is a natural goal for optimization and control. However, voltage variations along the circuit must stay within strict regulation bounds. Measured on a per unit basis, the voltage bounds become
\begin{eqnarray}
\forall j=0,\cdots,n:\quad 1-\epsilon\leq V_j^2\leq 1+\epsilon,
\label{Voltage}
\end{eqnarray}
where normally $\epsilon\approx 0.05$.

Combining all of the above,  we arrive at the following global Distflow optimization task \cite{89BWb,89BWa,89BWc}:
\begin{eqnarray}
\label{DistFlow}
\left.\min_{{\bm P},{\bm Q},{\bm V},{\bm q}^{(g)}} {\cal L}\right|_{\mbox{Eqs.~(\ref{Pj+1},\ref{Qj+1},\ref{Vj2},\ref{PV_constraint},\ref{Voltage})}},
\end{eqnarray}
for known impedances and given configuration of ${\bm p}^{(c)},{\bm p}^{(g)},{\bm q}^{(c)}$. The general DistFlow problem is not convex and may have multiple solutions.  We reduce the computational burden by considering the simplified DistFlow problem based on the DC approximation, i.e., LinDistFlow\cite{89BWb,89BWa,89BWc}
\begin{eqnarray}
\label{LinDistFlow}
&&\min_{{\bm P},{\bm Q},{\bm V},{\bm q}^{(g)}} \sum_{j=0}^{n-1} r_j\frac{P_j^2+Q_j^2}{V_j^2},\\
&&\mbox{s.t.}\quad\begin{array}{c}
P_{j+1}= P_j-p_{j+1}^{(c)}+p_{j+1}^{(g)},\\
Q_{j+1}=Q_j-q_{j+1}^{(c)}+q_{j+1}^{(g)},\\
V_{j+1}^2=V_j^2-2(r_jP_j+x_jQ_j),\\
\mbox{Eqs.~(\ref{PV_constraint},\ref{Voltage})}.
\end{array}
\nonumber
\end{eqnarray}
The formulation is a convex quadratic problem (as the quadratic cost function is convex and all the constraints are linear) with unique solution that can be computed efficiently.  We will argue in Section \ref{sec:rural} that this DC-based approximation is well justified for our example of a rural distribution circuit.

\section{Description of the prototypical rural distribution circuit}
\label{sec:rural}

To demonstrate the technique, we consider a model distribution circuit based loosely on one of 24 prototypical distribution circuits described in a taxonomy of distribution circuits.\cite{08SCCPET}  Our model represents a sparsely-loaded rural distribution circuit with a nominal line-to-neutral voltage of 7.2 kV.  The line impedance $(0.33+0.38i) \Omega /km$ is constant and based on typical conductor types, sizes and spacings.\cite{08SCCPET}

For this initial study, we consider 100 load nodes separated by distances uniformly distributed in the range $200$ to $300$ meters.  The real power consumed at each node ($p_j^{(c)}$) is selected from a uniform distribution between 0 and 4 kW, and the reactive power consumed ($q_j^{(c)}$) is selected from uniform distribution between $0.2 p_j^{(c)}$ and $0.3 p_j^{(c)}$ yielding power factors for each load distributed between approximately 0.96 to 0.98. At nodes with PV-generation capability, the real power generated ($p_j^{(g)}$) was always 1 kW reflecting a situation where the solar insolation is constant over the distribution circuit and all PV systems are the same size.  We explore different PV penetration levels by randomly assigning a variable fraction of the nodes to have PV generation.  At those nodes, we assign the same apparent power capacity $s$ and the reactive power generated (or consumed) is bounded by Eq.~\ref{PV_constraint}.  We explore the effects of different amounts of inverter excess apparent power capacity by varying $s$.

In our model the characteristic values of the nonlinear terms $\propto (P_j^2 +Q_j^2)/V_j^2$ in Eqs.~(\ref{Pj+1},\ref{Qj+1}) are about $10^{4}$ times smaller compared to the linear terms $P_j,Q_j$, so modeling based on the LinDistFlow DC approximation (\ref{LinDistFlow}) produces results almost indistinguishable from the exact AC model (\ref{Pj+1},\ref{Qj+1}).

\section{Simulations: Results and Discussions}
\label{sec:simulations}

For a given $s$ and PV-penetration factor $r$, we considered many different realizations of the prototype circuit.  Each realization consisted of different $p_j^{(c)}$, $q_j^{(c)}$, and distances between adjacent nodes.  The realizations are generated by drawing from the distributions described in the previous section.  By considering many samples for each combination of $s$ and $r$, we observe that the important qualitative and quantitative results described below are sufficiently robust.  Therefore, in this publication, we present results for a typical sample and leave the detailed discussion of statistics for further and more formal exploration.  Each realization serves as the entry point for the DistFlow or LinDistFlow optimization algorithm.  However, as argued above, the nonlinear terms in the AC power flow are quite small and LinDistFlow provides an accurate solution.  

\begin{figure}
\centering \caption{
Energy saved (in percent of the energy lost if all $q^{(g)}$ are set to zero) as a function of $s$. Here, $s$ is reported in kW, and $s=1.1$ implies that the apparent power capacity of the inverters is $10\%$ higher than the maximum real power output of the PV system.  The different curves correspond to configurations with different levels of PV penetration $r$.} \includegraphics[width=0.48\textwidth]{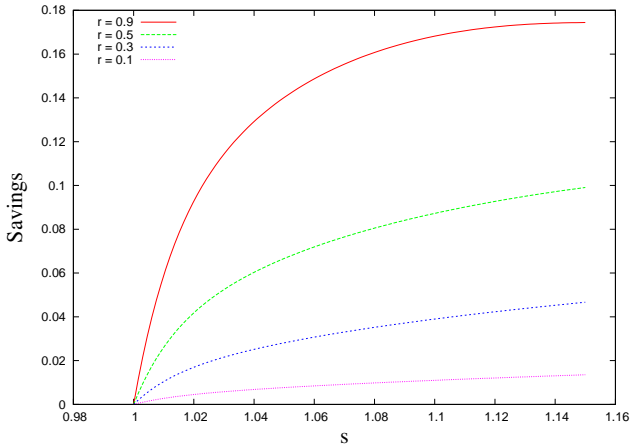}
\label{fig:loss_vs_sa}
\end{figure}

In our preliminary study of this problem, we investigate the reduction in losses first as a function of $s$ for several different values of $r$.  Second we investigate the reduction in loss as a function of $r$ for a single value of $s$.  We compare this result to that obtained via a simple distributed control approach.  Finally, we present the improvement in power quality (i.e. voltage regulation) obtained by using a centralized algorithm to dispatch $q_j^{(g)}$ from the nodes with PV capability.

Figure~\ref{fig:loss_vs_sa} shows that the energy savings increase monotonically with the inverter apparent power capacity $s$.  For a wide range of PV penetration, the majority of the energy savings occurs by $s=1.1$, i.e. at an apparent power capacity of only $10\%$ higher than $\max p_j^{(g)}$.  Additional savings are possible beyond $s=1.1$, however, higher values of $s$ will increase the cost of PV systems.  Determining the optimal value of $s$ must also include economic factors and is beyond the scope of this work.  The energy savings and the value of $s$ at which it saturates will depend on the reactive power consumed by the loads, however, we have not explored this dependence in this preliminary study.  Also, we have only considered the case when all the $p_j^{(g)}$ are at their maximum value.  We note that the saturation observed in Fig.~\ref{fig:loss_vs_sa} will occur faster in more realistic models of renewable generation where the PV systems will not always operate at their maximum capacity.

 \begin{figure}
\centering \caption{
Energy saved
(in percentage of the total energy loss observed when all $q^{(g)}$ are set to zero)
as a function of PV penetration measured in terms of the fraction  $r$ of nodes that can inject reactive power. Here, we set $s=1.1$.  The solid-red curve shows results of global optimization.  The dashed-green curve corresponds to the case of local control where the injected reactive power is equal to either the reactive power consumed or to $\sqrt{s_j^2-(p_j^{(g)})^2}$ if the consumed power is higher than the respective bound on $q^{(g)}$ given in Eq.~\ref{PV_constraint}} \includegraphics[width=0.48\textwidth]{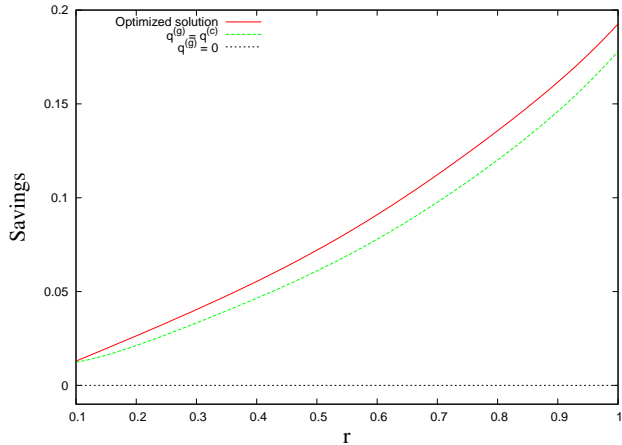}
\label{fig:loss_vs_r}
\end{figure}

Figure~\ref{fig:loss_vs_r} shows how the energy savings depend on PV penetration for $s=1.1$, i.e. a vertical slice through Fig~\ref{fig:loss_vs_sa} at $s=1.1$.  First we note that the energy savings approaches $20\%$ at very high PV penetrations suggesting that, if reactive power dispatch can be achieved, high PV penetration is in fact beneficial.  Again, we note that the quantitative results presented in this manuscript depend on our assumptions of the reactive power consumed by the loads.  However, we do not expect the qualitative behavior to change as the load assumptions are changed.  Figure~\ref{fig:loss_vs_r} also compares a global optimization (red) with a local control scheme\cite{08LB} requiring no communication at all.   In this particular local scheme, $q^{(g)}$ is set to minimize the local reactive power consumption $q_j$ subject to the constraints of Eq.~\ref{PV_constraint}. For the prototypical circuit and load assumptions used in this preliminary study, this naive control strategy surprisingly achieves about $95\%$ of the maximum possible savings.  The efficacy of this local scheme should be explored for many different circuit configurations and loads.  Note that both strategies should be compared with the ``do nothing", $q_j^{(g)}=0$ case shown in dashed black in Fig.~\ref{fig:loss_vs_r}.

\begin{figure}
\centering \caption{
Voltage variation and $q_j^{(g)}$ for a particular circuit realization with $r=90\%$ and $s=2.0$.   The voltage drop (normalized with respect to initial voltage $V_0$) is less if one is allowed to inject reactive power into the circuit.} \includegraphics[width=0.48\textwidth]{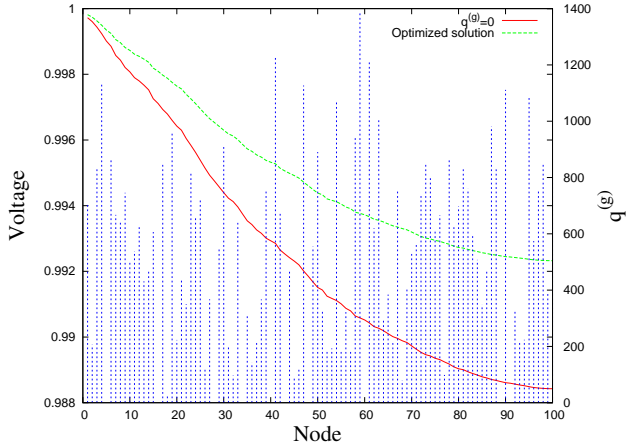}
\label{fig:land}
\end{figure}

Finally, we demonstrate the effect on power quality (i.e. voltage regulation) of dispatching reactive power from a relatively high penetration of PV inverters.  Figure~\ref{fig:land} shows the voltage profile along the circuit with (green) and without (red) reactive power dispatch.  Here, the reactive power is dispatched using the global scheme described above.  In addition to significant reduction in losses, the controlled dispatch of reactive power increases the power quality by significantly reducing the voltage drop along the circuit.  Although not shown in the Figure, we note that the global control of reactive power performs better than the local control because it allows for better suppression of reactive power flow $Q_j$ over the whole line, not just $q_j$.  Figure~\ref{fig:land} also shows the individual $q_j^{(g)}$ at each node along the line (vertical blue dashed lines, each standing for a generator with some lines missing because $r \leq 1.0$).

\section{Conclusions and Path Forward}
\label{sec:conclusions}
Our preliminary study is meant to elucidate the qualitative behavior of a distribution circuit with various levels of PV system penetration, excess inverter apparent power capacity, and global versus local inverter-reactive power dispatch.  We find that power dissipation is reduced and power quality is increased as  both PV penetration and excess apparent power capacity are increased.  For the prototypical circuit and loads considered in this study, we find that the reduction in power dissipation plateaus at a relatively low value of excess apparent power and that a local control scheme performs nearly as well as a global solution.

Perhaps more importantly, this study helps to formulate future research directions to
improve the optimization and control algorithms. We envision the following future modifications to the basic scheme discussed in the manuscript:
\begin{itemize}
\item Different circuit configurations and load profiles must be explored to determine if the qualitative and quantitative results apply to a wide variety of cases.  The taxonomy of distribution circuits described in \cite{08SCCPET} provides a guide to selecting additional case studies.
\item The global optimality of the algorithm can be traded for local optimization. The naive local scheme discussed in the manuscript, which requires making only local decisions, is already significantly better than no compensation at all.  However, there are all other intermediate strategies which may close the gap between the fully local (simple but strictly suboptimal) and fully global (optimal but complex and requiring significant communication and coordination) strategies.  We anticipate that a linear-scaling algorithm capable of efficiently solving the single-branch DistFlow optimization problem may be available.  We envision this type of optimal algorithm to have a form of dynamical programming \cite{bellman} with the global optimization replaced by a sequence of local optimizations advanced sequentially from the end of the circuit towards the entry-point (and/or in the reversed direction). These techniques could also be extended beyond the single-branch circuit to the case of a circuit with multiple branches, provided that the circuit remains a tree (no loops).

\item Accounting for loops in the graph constitutes an interesting algorithmic challenge that may have applications for highly meshed distribution circuits typical of urban systems.  In this case, solving the DistFlow equations will likely be of exponential complexity (in the number of nodes), however the LinDistFlow equations are still expected to be polynomially-tractable, i.e. reducible to a linear program. Normally, linear programming is not distributed, however in some special cases, e.g. minimum-cost network flow \cite{10GSW}, it allows a distributed implementation via Belief Propagation (BP) algorithms. Developing a distributed BP-algorithm for LinDistFlow optimization over a loopy graph is thus another challenging task for future exploration. The distributed nature of the algorithm will naturally enable the ability to carry out the global optimization via a sequence of local computations and communications limited only to neighboring loads on the circuit.

\item  Following the general arguments of \cite{08HOFO}, adding some switching capabilities can have a beneficial impact on distribution circuit losses for tree-like circuits and potentially for loopy circuits. Analysis of both algorithmic and phase transition aspects of switching for such a loopy circuit can be carried out in the spirit of \cite{09ZDC,09ZBC}.

\end{itemize}

\section*{Acknowledgment}

We are thankful to all the participants of the ``Optimization and Control for Smart Grids" LDRD DR project at Los Alamos
and Smart Grid Seminar Series at CNLS/LANL for multiple fruitful discussions and to Prof. Ross Baldick for attracting our attention to the DistFlow studies \cite{89BWa,89BWb,89BWc}. Research at LANL was carried out under the auspices of the National Nuclear Security Administration of the U.S. Department of Energy at Los Alamos National Laboratory under Contract No. DE C52-06NA25396. PS and MC acknowledges partial support of NMC via NSF collaborative grant CCF-0829945 on ``Harnessing Statistical Physics for Computing and Communications''.

\bibliographystyle{IEEEtran}
\bibliography{SmartGrid}

\begin{thebibliography}{10}
\providecommand{\url}[1]{#1}
\csname url@samestyle\endcsname
\providecommand{\newblock}{\relax}
\providecommand{\bibinfo}[2]{#2}
\providecommand{\BIBentrySTDinterwordspacing}{\spaceskip=0pt\relax}
\providecommand{\BIBentryALTinterwordstretchfactor}{4}
\providecommand{\BIBentryALTinterwordspacing}{\spaceskip=\fontdimen2\font plus
\BIBentryALTinterwordstretchfactor\fontdimen3\font minus
  \fontdimen4\font\relax}
\providecommand{\BIBforeignlanguage}[2]{{%
\expandafter\ifx\csname l@#1\endcsname\relax
\typeout{** WARNING: IEEEtran.bst: No hyphenation pattern has been}%
\typeout{** loaded for the language `#1'. Using the pattern for}%
\typeout{** the default language instead.}%
\else
\language=\csname l@#1\endcsname
\fi
#2}}
\providecommand{\BIBdecl}{\relax}
\BIBdecl

\bibitem{89BWb}
M.~Baran and F.~Wu, ``Optimal capacitor placement on radial distribution
  systems,'' \emph{Power Delivery, IEEE Transactions on}, vol.~4, no.~1, pp.
  725--734, Jan 1989.

\bibitem{89BWa}
------, ``Optimal sizing of capacitors placed on a radial distribution
  system,'' \emph{Power Delivery, IEEE Transactions on}, vol.~4, no.~1, pp.
  735--743, Jan 1989.

\bibitem{90BW}
R.~Baldick and F.~Wu, ``Efficient integer optimization algorithms for optimal
  coordination of capacitors and regulators,'' \emph{Power Systems, IEEE
  Transactions on}, vol.~5, no.~3, pp. 805--812, Aug 1990.

\bibitem{08BMOCKS}
\BIBentryALTinterwordspacing
M.~McGranaghan, T.~Ortmeyer, D.~Crudele, T.~Key, Smith, and J.Baker, ``Advanced
  grid planning and operation,'' NREL/SR-581-42294, Tech. Rep., 2008. [Online].
  Available:
  \url{http://www1.eere.energy.gov/solar/pdfs/advanced_grid_planning_operation%
s.pdf}
\BIBentrySTDinterwordspacing

\bibitem{90GKFE}
E.~Gulachenski, E.~J. Kern, W.~Feero, and A.~Emanuel, ``Photovoltaic generation
  effects on distribution feeders, volume 1: Description of the gardner,
  massachusetts, twenty-first century pv community and research program,'' EPRI
  report EL-6754, Tech. Rep., 1990.

\bibitem{08WNRN}
\BIBentryALTinterwordspacing
C.~Whitaker, J.~Newmiller, M.~Ropp, and B.~Norris, ``Distributed photovoltaic
  systems design and technology requirements,'' Sandia/SAND2008-0946 P, Tech.
  Rep., 2008. [Online]. Available:
  \url{http://www1.eere.energy.gov/solar/pdfs/distributed_pv_system_design.pdf}
\BIBentrySTDinterwordspacing

\bibitem{08LB}
\BIBentryALTinterwordspacing
E.~Liu and J.~Bebic, ``Distribution system voltage performance analysis for
  high-penetration photovoltaics,'' NREL/SR-581-42298, Tech. Rep., 2008.
  [Online]. Available: \url{http://www1.eere.energy.gov/solar/pdfs/42298.pdf}
\BIBentrySTDinterwordspacing

\bibitem{1547}
\BIBentryALTinterwordspacing
``{IEEE 1547 Standard for Interconnecting Distributed Resources with Electric
  Power Systems}.'' [Online]. Available:
  \url{http://grouper.ieee.org/groups/scc21/1547/1547_index.html}
\BIBentrySTDinterwordspacing

\bibitem{89BWc}
M.~Baran and F.~Wu, ``Network reconfiguration in distribution systems for loss
  reduction and load balancing,'' \emph{Power Delivery, IEEE Transactions on},
  vol.~4, no.~2, pp. 1401--1407, Apr 1989.

\bibitem{08SCCPET}
\BIBentryALTinterwordspacing
K.~Schneider, Y.~Chen, D.~Chassin, R.~Pratt, D.~Engel, and S.~Thompson,
  ``Modern grid initiative-distribution taxonomy final report,'' Tech. Rep.,
  2008. [Online]. Available:
  \url{http://www.gridlabd.org/models/feeders/taxonomy_of_prototypical_feeders%
.pdf}
\BIBentrySTDinterwordspacing

\bibitem{bellman}
R.~Bellman, \emph{Dynamic Programming}.\hskip 1em plus 0.5em minus 0.4em\relax
  Princeton University Press, Princeton, NJ, 1957.

\bibitem{10GSW}
D.~Gamarnik, D.~Shah, and Y.~Wei, ``Belief propagation for min-cost flow:
  convergence and correctness,'' in \emph{ACM-SIAM Symposium on Discrete
  Algorithms (SODA10)}, 2010.

\bibitem{08HOFO}
K.~Hedman, R.~O'Neill, E.~Fisher, and S.~Oren, ``Optimal transmission
  switching—sensitivity analysis and extensions,'' \emph{Power Systems, IEEE
  Transactions on}, vol.~23, no.~3, pp. 1469--1479, Aug. 2008.

\bibitem{09ZDC}
\BIBentryALTinterwordspacing
L.~Zdeborova, A.~Decelle, and M.~Chertkov, ``Message passing for optimization
  and control of power grid: Model of distribution system with redundancy,''
  \emph{Physical Review E}, vol.~80, 2009. [Online]. Available:
  \url{http://arxiv.org/abs/0904.0477}
\BIBentrySTDinterwordspacing

\bibitem{09ZBC}
\BIBentryALTinterwordspacing
L.~Zdeborova, S.~Backhaus, and M.~Chertkov, ``Message passing for integrating
  and assessing renewable generation in a redundant power grid,'' 2009.
  [Online]. Available: \url{http://arxiv.org/abs/0909.2358}
\BIBentrySTDinterwordspacing

\end{thebibliography}

\end{document}